\documentclass[11pt]{article}%
\usepackage[singlespacing]{setspace}
\usepackage{amsfonts}
\usepackage{graphicx}
\usepackage{amsmath}
\usepackage{amssymb}
\usepackage{layout}%
\providecommand{\U}[1]{\protect\rule{.1in}{.1in}}
\newtheorem{theorem}{Theorem}

\newtheorem{lemma}[theorem]{Lemma}

\newtheorem{proposition}[theorem]{Proposition}
\newtheorem{remark}[theorem]{Remark}

\newcommand\E{{\mathbb E}}
\newcommand\p{{\mathbb P}}

\singlespace
\begin{document}

\title{On the central limit theorem for stationary random fields under  ${\mathbb L}^1$-projective condition }
\author{Han-Mai Lin{\thanks{Universit\'{e} Gustave Eiffel, LAMA and CNRS
UMR 8050. Email: han-mai.lin@univ-eiffel.fr}}, Florence Merlev\`{e}de{\thanks{Universit\'{e} Gustave Eiffel, LAMA and CNRS
UMR 8050. Email: florence.merlevede@univ-eiffel.fr}} and Dalibor Voln\'y\thanks{Universit\'e de Rouen, LMRS and CNRS UMR 6085. 
Email: dalibor.volny@univ-rouen.fr}}
\maketitle

\abstract{The first aim of this paper is to wonder to what extent we can generalize the central limit theorem of Gordin  \cite{Go73} under the so-called ${\mathbb L}^1$-projective criteria to  ergodic stationary random fields when completely commuting filtrations are considered. Surprisingly it appears that this result cannot be extended to its full generality and that an additional condition is needed.}

\medskip

\textit{Keywords:} Central limit theorem; stationary random fields; orthomartingale; coboundary decomposition;  projective criteria.

\smallskip

\textit{MSC2020:} 60F05; 60G60.
\bigskip

\section{Introduction and main results}

Let $(\Omega,{\mathcal A},  \p)$ be a probability space, and $T:\Omega
\mapsto \Omega$ be
 an {\it ergodic} bijective bimeasurable transformation preserving the probability $\p$. Let ${\cal F}_0 $ be a
sub-$\sigma$-algebra of ${\cal A}$ satisfying ${\mathcal F}_0
\subseteq T^{-1 }({\cal F}_0)$ and $f$ be a ${\mathbb L}^1(\p)$ real-valued centered random variable adapted to ${\cal F}_0 $. By $U$ we denote 
the operator $U$ : $f\mapsto f\circ T$. The notation $I$ will denote the identity operator. Define then 
the stationary sequence $(f_i)_{ i \in {\mathbb Z}}$ by $f_i=f
\circ T^i=U^i f$,  its associated stationary filtration $({\cal F}_i)_{ i \in {\mathbb Z}}$ by ${\mathcal F}_i={\mathcal F}_0
\circ T^{-i}$ and let $
   S_n(f) = \sum_{i=0}^{n-1} U^i f   $. 
   
 The following theorem is essentially due to Gordin \cite{Go73} and gives sufficient conditions for $(U^i f)_{ i \in {\mathbb Z}}$ to satisfy the central limit theorem. 
\begin{theorem}[Gordin] \label{Gordin} Assume that the series
\begin{equation} \label{L1-cond}
\sum_{i \geq 0} \E ( U^i f | {\mathcal F}_0) \text{ converges in ${\mathbb L}^1(\p)$}
\end{equation}
and
\begin{equation} \label{condliminf}
\liminf_{n \rightarrow \infty} \frac{\E (|S_n(f) |)}{ \sqrt{n}} < \infty \, .
\end{equation}
Then $n^{-1/2} S_n(f)$ converges in  distribution to a centered normal variable (that can be degenerate). 
\end{theorem}
The proof of this result is based on the following coboundary martingale decomposition (see \cite{V13} for more details concerning necessary and sufficient conditions for the existence of such a decomposition): Under \eqref{L1-cond}, 
\begin{equation} \label{co-mart}
f = m + (I-U)g 
\end{equation}
where $m$ and $g$ are in ${\mathbb L}^1 ( \p)$ and $(U^i m)_{i\geq 0}$ is a stationary sequence of martingale differences, and on the following theorem whose 
complete proof can  be found in Esseen and Janson \cite{EJ-85}. 

\begin{theorem}[Esseen-Janson]  \label{EJth} If $(U^i m)_{i\geq 0}$ is a stationary and ergodic sequence of martingale differences in ${\mathbb L}^1 (\p)$ satisfying 
\begin{equation} \label{condliminfMart}
\liminf_{n \rightarrow \infty} \frac{\E (| \sum_{i=0}^{n-1} U^i m |)}{ \sqrt{n}} < \infty \, ,
\end{equation}
then $m \in {\mathbb L}^2 ( \p)$. 
\end{theorem}
Clearly, using the coboundary martingale decomposition \eqref{co-mart},  \eqref{condliminf} implies \eqref{condliminfMart}. 

The aim of this paper is to prove that Theorem \ref{Gordin} can be extended to random fields when the underlying filtrations are completely commuting. To fix the idea, let us first state the result in case of multidimensional index of dimension $d=2$ (the general case will be stated in Section \ref{sectiongeneral}). Then, in complement to the previous notation,  let $S$ be an ergodic  bimeasurable and measure preserving bijection of $\Omega$.   By $V$
we denote the operator $V$ : $f\mapsto f\circ S$. 

{\it In what follows we shall assume that the ergodic transformations $T$ and $S$ are commuting}. Note that  $T_{i, j} = T^iS^j$ is an ergodic ${\mathbb Z}^2$ action on $(\Omega, \mathcal A, \p)$. Let ${\mathcal F}_{0,0}$ be a sub-sigma field of ${\mathcal A}$ and for all $(i,j) \in {\mathbb Z}^2$ define 
\[
{\mathcal F}_{i,j} = T^{-i} S^{-j} ({\mathcal F}_{0,0} ) \, .
\]
Suppose that the filtration $({\mathcal F}_{i,j})_{ (i,j \in {\mathbb Z}^2} $  is increasing in $i$ for every $j$ fixed and increasing in $j$ for every $i$ fixed, and  is 
 {\it  completely commuting} in the sense that, for any integrable  $f$,
\[
{\mathbb E}( {\mathbb E}(f| {\mathcal F}_{i,j} ) | {\mathcal F}_{u,v} ) = {\mathbb E}(f| {\mathcal F}_{i\wedge u,j\wedge v} )  \, .
\]

\smallskip

In view of giving an extension of  Theorem \ref{Gordinext} for random fields indexed by the lattice ${\mathbb Z}^2$, the first tool is a suitable coboundary orthomartingale decomposition:  Let $f$ be a $\mathcal F_{0,0}$-measurable centered ${\mathbb L}^1(\p)$ function. According to Voln\'y \cite{V18},   the condition:
\begin{equation} \label{L1-cond-field}
\text{ the series }\sum_{i ,j\geq 0} \E ( U^i V^j f | {\mathcal F}_{0,0}) \text{ converges in ${\mathbb L}^1(\p)$}
\end{equation}
implies  the existence of the following decomposition: 
\begin{equation} \label{co-martext}
  f = m + (I - U)g_1 + (I - V)g_2  + (I-U)(I-V)g_3 \, , 
\end{equation}
where $m, g_1, g_2, g_3 \in {\mathbb L}^1(\p)$, $(U^iV^jm)$ is a stationary field of orthomartingale differences, $(V^jg_1)_j$ is a stationary martingale differences sequence 
with respect to the filtration $(\mathcal F_{\infty,j})_j$, and $(U^ig_2)_i$ is a stationary martingale differences sequence 
with respect to the filtration $(\mathcal F_{i,\infty})_i$.  
 To fix the ideas, setting $\E_{a,b} (\cdot) =  \E ( \cdot | {\mathcal F}_{a,b})$,  we have
\[
m = \sum_{i,j \geq 0} \big (   \E_{0,0} ( U^i V^j f ) - \E_{-1,0} ( U^i V^j f )  -  \E_{0,-1} ( U^i V^j f )  +  \E_{-1,-1} ( U^i V^j f )  \big ) \, , 
\]
\[
g_1 = \sum_{i,j \geq 0} \big (  \E_{-1,0} ( U^i V^j f )  -    \E_{-1,-1} ( U^i V^j f )  \big ) \, , \, g_2 = \sum_{i,j \geq 0} \big (  \E_{0,-1} ( U^i V^j f )  -    \E_{-1,-1} ( U^i V^j f )  \big ) \, , 
\]
and $g_3 = \sum_{i,j \geq 0}  \E_{-1,-1} ( U^i V^j f )  $.  Recall also that $(U^iV^jm)$  is said to be a orthomartingale differences field w.r.t. $ ({\mathcal F}_{i,j} )$ if 
\[
\E_{i-1,j} ( U^iV^jm ) = \E_{i,j-1} ( U^iV^jm )=\E_{i-1,j-1} ( U^iV^jm )= 0 \ \text{a.s.}
\]
Note that if  $f$ is additionally  assumed  to be regular in the sense that $f$ is $ {\mathcal F}_{\infty, \infty}$-measurable and $\E ( f | {\mathcal F}_{0,-\infty})= \E ( f | {\mathcal F}_{-\infty,0}) =0$ then,  it is proved in Voln\'y \cite{V18} that the converse is true, meaning that if $f$ satisfies the decomposition  \eqref{co-martext} then  \eqref{L1-cond-field} holds.  We also refer to \cite{EG} where the existence of the decomposition 
\eqref{co-martext} is proved under a reinforcement of  \eqref{L1-cond-field} (they assume that the series of the ${\mathbb L}^1$-norm is convergent). We also mention \cite[Theorem 2.2]{Gi18} where a necessary and sufficient condition for an orthomartingale-coboundary decomposition is established when all the underlying random elements are square integrable. 

\smallskip

Our first result is the following: 

\begin{theorem} \label{Gordinext} Let  $f$ be a $\mathcal F_{0,0}$-measurable centered ${\mathbb L}^1(\p)$ random variable.  Let $S_{n_1,n_2} (f) = \sum_{i=0}^{n_1}  \sum_{j=0}^{n_2}
U^i V^j f$. Assume that condition \eqref{L1-cond-field} is satisfied and that 
\begin{equation} \label{condliminf-field-1}
  \liminf_{n\to\infty}\frac1{\sqrt{n}}  \Big \Vert \sum_{i=0}^{n-1} 
  U^if \Big \Vert_1 <\infty \, ,  \,    \liminf_{N\to\infty}\frac1{\sqrt{N}} \Big \Vert \sum_{j=0}^{N-1} 
  V^j f \Big \Vert_1 <\infty \, 
  \end{equation}
  and
\begin{equation} \label{condliminf-field-2}
\liminf_{ n_1 \rightarrow \infty } \liminf_{ n_2 \rightarrow \infty } \frac{\E (|S_{n_1,n_2} (f)  |)}{ \sqrt{n_1 n_2}} < \infty  \, , \,  \liminf_{ n_2 \rightarrow \infty } \liminf_{ n_1 \rightarrow \infty } \frac{\E (|S_{n_1,n_2} (f)  |)}{ \sqrt{n_1 n_2}} < \infty  \, .
\end{equation}
Then the random variables $m$,   $(I - U)g_1$ and $  (I - V)g_2$ defined in \eqref{co-martext} are in ${\mathbb L}^2(\p)$.
\end{theorem}
Compared to the case of random sequences a natural question  is then to wonder if condition  \eqref{L1-cond-field} together with conditions \eqref{condliminf-field-1} and \eqref{condliminf-field-2} are sufficient to ensure that, when $\min (n_1 , n_2)  \rightarrow \infty$, the limiting distribution behavior of $(n_1n_2)^{-1/2} S_{n_1,n_2}(f)$ is the same as that of the orthomartingale part $(n_1n_2)^{-1/2} S_{n_1,n_2}(m)$. In other terms one can wonder if assuming  the conditions of Theorem \ref{Gordinext}  is enough to ensure that the coboundaries' behavior, i.e. $(n_1n_2)^{-1/2} ( S_{n_1,n_2} (f) -  S_{n_1,n_2}(m) )$ is negligible for the convergence in distribution. Surprisingly the answer to this question is negative as shown by the next counterexample. 

 \begin{theorem} \label{GordinextCE}  There exist a probability space $(\Omega, {\mathcal A}, \mu)$, a function  $g \in {\mathbb L}^1(\mu) $, measurable with respect to a $\sigma$-algebra ${\mathcal F}_{0,0} \subset {\mathcal A}$ and bijective bimeasurable ergodic transformations $T$ and $S$ such that $f = (I-U)g $ is in ${\mathbb {\mathbb L}^2}(\mu) $, satisfies the conditions   \eqref{L1-cond-field}, \eqref{condliminf-field-1} and  \eqref{condliminf-field-2}  but such that $(n_1n_2)^{-1/2} S_{n_1,n_2}(f)$ does not converge in  distribution to zero as $\min(n_1,n_2) \to \infty$. 
\end{theorem}
This result proves a drastically different behavior for the case of random fields with dimension $d \geq 2$ compared to the case of random sequences $d=1$ for which the coboundary is negligible for the convergence in distribution as soon as \eqref{L1-cond} is assumed. 
\begin{remark}
Modifying the selection of the sequences  $ (n_k) $ and $ (m_k)$ used in the construction of the counterexample of Theorem \ref{GordinextCE}, we infer that one can construct the function $g$ in such a way that not only it satisfies the conditions and the conclusion of  Theorem \ref{GordinextCE}  but also the following conditions:  for $p \in [1,2)$, $g \in {\mathbb L}^p(\mu)  $ and the series $\sum_{i ,j\geq 0} \E ( U^i V^j f | {\mathcal F}_{0,0})$ converges in ${\mathbb L}^p(\mu) $ (it suffices to take for instance $n_k= [2^{k/2} ] $ and $m_k \sim (n_k /k)^{p/(2-p)}$). 
\end{remark}

However, reinforcing the conditions of Theorem \ref{Gordinext}, we can prove the following CLT.

 \begin{theorem} \label{Gordinextpositive}  In addition to the conditions of Theorem  \ref{Gordinext}, assume that 
 \begin{equation} \label{condlimSomme}
\lim_{\min(n_1,n_2) \rightarrow \infty } \frac{\E (|S_{n_1,n_2} (f)  |)}{ \sqrt{n_1 n_2}}  \text{ exists.}  
\end{equation}
Then, as $\min(n_1,n_2) \to \infty$, $(n_1n_2)^{-1/2} S_{n_1,n_2}(f)$ converges in  distribution to a centered normal variable (that can be degenerate).
\end{theorem}
\begin{remark} Assume that $f$ satisfies the coboundary orthomartingale decomposition \eqref{L1-cond-field} with the following additional conditions: 
$(U^iV^jm)$ is a stationary field of ${\mathbb L}^2(\p)$ orthomartingale differences, $(V^j (I-U)g_1)_j$ is a ${\mathbb L}^2(\p)$ stationary martingale differences sequence 
with respect to the filtration $(\mathcal F_{\infty,j})_j$, and $(U^i(I-V)g_2)_i$ is a ${\mathbb L}^2(\p)$ stationary martingale differences sequence 
with respect to the filtration $(\mathcal F_{i,\infty})_i$ (Theorem \ref{Gordinext}  gives sufficient conditions ensuring such a decomposition). Then, from the proof of Theorem \ref{Gordinextpositive}, we infer that condition \eqref{condlimSomme} is equivalent to the two following conditions:
\begin{align*}
 \limsup_{n\to\infty} \limsup_{k\to\infty} \frac1{\sqrt{nk}} \|S_{n,k}(f)\|_1 &  \leq \lim_{n\to\infty} \frac1n \|S_{n,n}(f)\|_1, \\
  \limsup_{k\to\infty} \limsup_{n\to\infty} \frac1{\sqrt{nk}} \|S_{n,k}(f)\|_1&   \leq \lim_{n\to\infty} \frac1n \|S_{n,n}(f)\|_1.
\end{align*} 
Note that the existence of the limit $ \lim_{n\to\infty} \frac1n \|S_{n,n}(f)\|_1$ has been mentioned in the proof of Theorem \ref{Gordinextpositive}.
\end{remark}
It is noteworthy to indicate that $f$ does not need to be in ${\mathbb L}^2$ but only  in ${\mathbb L}^1$  to apply Theorem \ref{Gordinextpositive} (see Example \ref{ex1} given below).  Theorem \ref{Gordinextpositive}  then gives alternative projective conditions compared to those required in \cite[Th. 5.1]{VW14} or in \cite[Th. 1]{PZ17} for the central limit theorem under the normalization $\sqrt{n_1n_2}$ to hold.  Note that the proofs of the two above mentioned  results  are also based on an orthomartingale approximation.  We refer also to \cite{WV13} where the notion of orthomartingales and completely commuting filtrations have been previously used in the particular case of functions of iid random fields.  Let us also indicate that when filtrations in the lexicographic order rather than completely  commuting filtrations are considered, \cite[Th. 1]{D98} provides a projective type condition in the spirit of the ${\mathbb L}^1$-projective condition  \eqref{L1-cond-field} (but still requiring $f$ to be in ${\mathbb L}^2$) for the normalized partial sums associated with a stationary random field to satisfy the central limit theorem. His proof is based on the so-called Lindeberg method.

\section{Examples}

\subsection{An  example of application when $f$ is in ${\mathbb L}^1$ but not in ${\mathbb L}^2$} \label{ex1}

 For $k \in {\Bbb N}^*$ and $i,j\in \Bbb Z$, let $e_{k,i,j}$ be mutually independent zero mean random variables with
$Ue_{k,i,j} = e_{k,i+1,j}$, $V e_{k,i,j} = e_{k,i,j+1}$. Let ${\mathcal F}_{a,b} = \sigma ( e_{k,i,j} , k \in {\mathbb Z},  i \leq a , j \leq b )$. We denote
$e_k = e_{k,0,0}$. Assume that  for any $i,j\in \Bbb Z^2$, ${\mathcal L}(e_{k,i,j} ) = {\mathcal L}(e_{k} ) $ and that $e_k$ takes value $v_k$ with probability $p_k$, $-v_k$ with probability $p_k$ and $0$ with probability $1-2p_k$. It follows that 
\[ 
\Vert e_k \Vert_1 = 2 v_k p_k \text{ and }  \Vert e_k \Vert_2^2 = 2 v^2_k p_k \, .
\]
We choose $(v_k)_{k \geq 1} $ and $(p_k)_{k \geq 1}$ as follows
\[
v_k = k^2  (\log (k+1) )^2 \text{ and }   p_k =  \frac{1}{ 2 k^2 (\log (k+1) )^4} \, .
\]
For any $k \geq 1 $ and $i \geq 0$, let $ \displaystyle 
a_{k,i} = \frac{1}{(k+i)^2} $, and define 
\begin{align*}
 &  g_1 = \sum_{k\geq 1}  \sum_{i \geq 0} a_{k,i} U^{-k-i} e_k, \quad g_2 =  \sum_{k\geq 1}  \sum_{j \geq 0} a_{k,j} V^{-k-j} e_k, \quad
  g_3 =  \sum_{k=1}^\infty \frac1k e_k , \\ 
&  m = \sum_{k=1}^\infty \frac{1}{k^2} U^{-k} V^{-k}  e_k \,  \text{ and }  \,   f = m +  (I-U)g_1 + (I-V)g_2 + (I-U)(I-V)g_3.
\end{align*}
It follows that 
\[
\Vert   g_1 \Vert_1  \leq  \sum_{k\geq 1}  \sum_{i \geq k}  \frac{1}{ i^2}\Vert e_k  \Vert_1 \leq 2 \sum_{k\geq 1}    \frac{1}{ k }\Vert e_k  \Vert_1  \leq   \sum_{k\geq 1}^\infty \frac{2}{k   (\log (k+1))^2}   < \infty \, .
\]
But, by independence of the r.v.'s $e_{k,i,j}$, 
\[
\Vert   g_1 \Vert^2_2  =  \sum_{k\geq 1 }   \sum_{i \geq 0} a^2_{k,i} \Vert e_k  \Vert_2^2  =     \sum_{k\geq 1 }   \sum_{i \geq k}  \frac{1}{i^4}\Vert e_k  \Vert_2^2  \geq  \sum_{k\geq 1} \frac{1 }{ 3 k}   = \infty \, .
\]
Similarly $\Vert   g_2 \Vert_1+ \Vert   g_3 \Vert_1 < \infty$, $\Vert   g_2 \Vert_2 = \infty $ and $\Vert   g_3 \Vert_2 =  \infty$. 
On another hand, for any positive integer $\ell$, 
\[
(I-U^{\ell})g_1 = -  \sum_{k\geq 1}  \sum_{i =0}^{\ell-1} a_{k,i} U^{-k-i+\ell} e_k + \sum_{k=1}^\infty   \sum_{i \geq 0}(a_{k,i} - a_{k,i+\ell }  ) U^{-k-i} e_k \, .
\]
Hence,  by independence of the r.v.'s $e_{k,i,j}$, 
\begin{align*}
\Vert (I-U^{\ell})g_1 \Vert_2^2 & =   \sum_{k\geq 1}  \sum_{i =0}^{\ell-1} a^2_{k,i}  \Vert  e_k  \Vert_2^2 + \sum_{k=1}^\infty   \sum_{i \geq 0}(a_{k,i} - a_{k,i+\ell }  )^2  \Vert  e_k  \Vert_2^2  \\
& =   \sum_{k\geq 1}  k^2 \sum_{i =0}^{\ell-1} a^2_{k,i}   + \sum_{k=1}^\infty  k^2  \sum_{i \geq 0} (a_{k,i} - a_{k,i+\ell }  )^2   \, .
\end{align*}
By simple algebra, we derive that there exists a positive constant $C$ such that for any positive integer $\ell$,
\begin{equation} \label{cobordmixte}
\Vert (I-U^{\ell})g_1 \Vert_2^2 \leq C  \log ( \ell +1 )  \, .
\end{equation}
In particular, we get $\Vert   (I-U) g_1 \Vert_2  < \infty$.  Similarly, we have $\Vert   (I-V) g_2 \Vert_2  < \infty$. 
To summarize, we have $g_1, g_2, g_3, (I-U)(I-V)g_3, f \in {\mathbb L}^1\setminus {\mathbb  L}^2$, $(I-U)g_1, (I-V)g_2 \in {\mathbb L}^2$.  On another hand $(V^j\big[(I-U)g_1\big])_j$ are martingale
differences as well as $(U^i\big[(I-V)g_2\big])_i$. \newline
For $f$, we clearly get   $\sum_{i,j \geq 0} \Vert {\mathbb E} ( U^i V^j f | {\mathcal F}_{0,0}  ) \Vert_1 < \infty $. In addition 
\[
\frac{\Vert \sum_{i=0}^{n-1} \sum_{j=0}^{N-1} U^i  V^j (I-U) g_1 \Vert_2^2 }{nN} = \frac{\Vert  (I-U^n) g_1 \Vert_2^2 }{n} \, .
\]
Taking into account \eqref{cobordmixte}, it follows that 
\[
 \lim_{n \rightarrow \infty }\frac{\Vert \sum_{i=0}^{n-1} \sum_{j=0}^{N-1} U^i  V^j (I-U) g_1 \Vert_2^2 }{nN} =0  \, .
\]
Similarly
\[
 \lim_{N \rightarrow \infty }\frac{\Vert \sum_{i=0}^{n-1} \sum_{j=0}^{N-1} U^i  V^j (I-V) g_2 \Vert_2^2 }{nN}  =0  \, .
\]
Moreover, as $ \max( n, N )  \rightarrow \infty$, the coboundary is negligible in ${\mathbb L}^1$, i.e.
 \[
 \lim_{\max( n, N )  \rightarrow \infty }\frac{\Vert \sum_{i=0}^{n-1} \sum_{j=0}^{N-1} U^i  V^j (I-U)(I-V) g_3 \Vert_1 }{ \sqrt{nN}}  =0  \, .
\]
It follows that
\[
 \lim_{\min(n,N) \rightarrow \infty }\frac{\Vert \sum_{i=0}^{n-1} \sum_{j=0}^{N-1} U^i  V^j (f-m)  \Vert_1 }{\sqrt{nN}}  =0  \, . 
\]
Since, $
 \lim_{\min(n,N) \rightarrow \infty }\frac{\Vert \sum_{i=0}^{n-1} \sum_{j=0}^{N-1} U^i  V^j m   \Vert_1 }{\sqrt{nN}} $ exists (it is equal to $\sqrt{\frac{2}{\pi}} \Vert m \Vert_2$), we can deduce that $
 \lim_{\min(n,N) \rightarrow \infty }\frac{\Vert \sum_{i=0}^{n-1} \sum_{j=0}^{N-1} U^i  V^j f   \Vert_1 }{\sqrt{nN}} $ exists. It follows that, all the conditions of Theorem \ref{Gordinextpositive}  are satisfied, and consequently, as $\min(n_1,n_2) \to \infty$, $(n_1n_2)^{-1/2} S_{n_1,n_2}(f)$ converges in  distribution to a centered normal variable. 
 
 \subsection{An example where $f$ does not satisfy Hannan's ${\mathbb L}^2$-condition} \label{ex2}  We exhibit an example where $f$ satisfies all the conditions of Theorem  \ref{Gordinextpositive} but not the Hannan's ${\mathbb L}^2$-condition required in  \cite[Th. 5.1]{VW14}.

We consider the random field $(e_{k,i,j})_{k,i,j}$ of mutually independent random variables as in  Example \ref{ex1} with the following choices of  $(v_k)_{k \geq 1}$ and $(p_k)_{k \geq 1} $. Let  $\alpha >4$. Then for any $k\geq 1$, define   
\begin{equation*}
		v_k = k^{\alpha}\text{ and }p_k = \dfrac{1}{2k^5\log(k+1)^2}\, .
	\end{equation*}
	Therefore
	\begin{equation*}
		\Vert e_k \Vert_1= \dfrac{k^{\alpha - 5}}{\log(k + 1)^2}  \text{ and } \Vert e_k \Vert_2^2 = \dfrac{k^{2\alpha - 5}}{\log(k + 1)^2} \, .
	\end{equation*}
For any $k \geq 1$ and $i,j \geq 0$, let 
\[
a_{k,i,j} = \frac{1}{(k+i+j)^{\alpha}} \, .
\]
Then, define 
\[
f = \sum_{k \geq 1} \sum_{u,v \geq 0} a_{k,u,v} U^{-u} V^{-v} e_k \, .
\]
$(U^i V^j f)_{i,j}$ is usually called a super linear random field. 
Let ${\mathcal F}_{a,b} = \sigma ( e_{k,i,j} , k \in {\mathbb Z},  i \leq a , j \leq b )$. 
Note that  $f$ is a $\mathcal F_{0,0}$-measurable centered and ${\mathbb L}^2(\p)$ random variable. In addition
\[
 \sum_{k \geq 1} \sum_{i,j \geq 0} \sum_{u \geq i+1}   \sum_{v \geq j+1}  |a_{k,u,v}  | \Vert e_k \Vert_1 \leq C   \sum_{k \geq 1} \frac{1}{k \log(k + 1)^2 }< \infty \, .
\]
Therefore  condition \eqref{L1-cond-field} is satisfied. So, the orthomartingale coboundary decomposition \eqref{co-martext} holds with 
\[
		m = \sum_{k \geq  1}\sum_{i ,j  \geq 0 } a_{k,i,j} e_k \, , 
		\]
		\[
		g_1 = \sum_{k \geq  1}\sum_{i ,j  \geq 0 }  \sum_{u \geq i+1}a_{k,u,j}  U^{i-u}e_{k}\, , \, 
		g_2 = \sum_{k \geq  1}\sum_{i ,j  \geq 0 }  \sum_{v \geq j+1}a_{k,i,v}  V^{j-v}e_{k}\, , 
		\]
		and
		\[
			g_3 = \sum_{k \geq  1}\sum_{i ,j  \geq 0 }  \sum_{u \geq i+1} \sum_{v \geq j+1}a_{k,u,v}  U^{i-u} V^{j-v}e_{k}\, .
	\]
One can easily check that $m \in {\mathbb L}^2(\p)$. Next, for any positive integer $\ell$, note that 
\[
(I-U^{\ell}) g_1 =  -  \sum_{k \geq  1}\sum_{i ,j  \geq 0 }  \sum_{u = 1}^{\ell +1}a_{k,u+i,j}  U^{-u}e_{k} + \sum_{k \geq  1}\sum_{i ,j  \geq 0 }  \sum_{u  \geq 1} ( a_{k,u +i,j}  - a_{k,u +i + \ell ,j}  )  U^{-u} e_{k} \, .
\]
By independence of the r.v.'s $e_{k,i,j}$,  it follows that 
\begin{align*}
\Vert (I-U^{\ell})g_1 \Vert_2^2 & \leq 2   \sum_{k\geq 1}    \sum_{u = 1}^{\ell +1}  \Big (  \sum_{i ,j  \geq 0 } a_{k,u+i,j} \Big )^2 \Vert  e_k  \Vert_2^2  \\ 
&  \quad \quad \quad   +   2 \sum_{k\geq 1}    \sum_{u \geq 1}  \Big (  \sum_{i,j \geq 0}  ( a_{k,u +i,j}  - a_{k,u +i + \ell ,j}  )\Big )^2 \Vert  e_k  \Vert_2^2    \, .
\end{align*}
Therefore,  there exists a positive constant $K_1$ (depending on $\alpha$) such that 
\begin{multline*}
\Vert (I-U^{\ell})g_1 \Vert_2^2   \\ \leq  K \Big (    \sum_{k\geq 1}   \sum_{u = 1}^{\ell +1}  \frac{1}{ (k +u )^{2 (\alpha -2)}}  \Vert  e_k  \Vert_2^2  + \ell^2   \sum_{k\geq 1}   \sum_{u \geq  1} \frac{1}{ (k +u + \ell )^{2 } (k +u  )^{2 (\alpha -2)}}  \Vert  e_k  \Vert_2^2  \Big )    \, .
\end{multline*}
By simple algebra, we then derive that there exists a positive constant $K_2$ (depending on $\alpha$) such that  such that for any positive integer $\ell$,
\begin{equation} \label{cobordmixte-2}
\Vert (I-U^{\ell})g_1 \Vert_2^2 \leq K_2 \frac{ \ell } {  \log ( \ell +1 )  } \, .
\end{equation}
Similarly, we get that there exists a positive constant $K_3$ (depending on $\alpha$) such that  such that for any positive integer $\ell$,
\begin{equation} \label{cobordmixte-3}
\Vert (I-V^{\ell})g_2 \Vert_2^2 \leq K_3 \frac{ \ell } {  \log ( \ell +1 )  } \, .
\end{equation}
From  \eqref{cobordmixte-2} and \eqref{cobordmixte-3}, it follows that 
\begin{equation} \label{cobordmixte-4}
 \lim_{n \rightarrow \infty }\frac{\Vert (I-U^{n})g_1 \Vert_2^2 +\Vert (I-V^{n})g_2 \Vert_2^2 }{n} =0  \, .
\end{equation}
Taking into account \eqref{cobordmixte-4} together with the fact that  $m \in {\mathbb L}^2(\p)$ and proceeding as in Example \ref{ex1}, one can verify that conditions  \eqref{condliminf-field-1},  \eqref{condliminf-field-2} and \eqref{condlimSomme} are satisfied.  Hence, Theorem  \ref{Gordinextpositive}   applies. Then, as $\min(n_1,n_2) \to \infty$, $(n_1n_2)^{-1/2} S_{n_1,n_2}(f)$ converges in  distribution to a centered normal variable. 

\medskip

On another hand, defining 
\[
P_{0,0} ( \cdot  ) = \E_{0,0} ( \cdot  )  -  \E_{-1,0} ( \cdot  )  -  \E_{0,-1} ( \cdot  )  +  \E_{-1,-1} ( \cdot  )  \, , 
\]
we get, for any $i,j \geq 0$,  $P_{0,0} ( U^i V^j f ) = \sum_{k \geq 1} a_{k,i,j} e_k$.  Hence, for any $i,j \geq 0$,
\[
\Vert  P_{0,0} ( U^i V^j f )  \Vert_2^2 =\sum_{k \geq 1} a^2_{k,i,j} \Vert  e_k  \Vert_2^2  \geq  \sum_{k \geq i+j +1}   \frac{k^{2 \alpha -5}}{(k+i+j)^{2 \alpha}  ( \log (k+1))^2}    \, ,
\]
implying that 
\[
\Vert  P_{0,0} ( U^i V^j f )  \Vert_2^2  \geq   \frac{C }{(i+j +1)^{ 4}  ( \log (i+j+2))^2}  \, .
\]
It follows that the series 
$
\sum_{i,j} \Vert  P_{0,0} ( U^i V^j f )  \Vert_2 $ diverges and then the Hannan's  ${\mathbb L}^2$-condition in the random fields setting does not hold, and \cite[Th. 5.1]{VW14} does not apply.  Note also that for this example, \cite[Th. 1]{D98} that involves filtrations in the lexicographic order, cannot be applied. 

\section{Proofs}

\subsection{Proof of Theorem \ref{Gordinext}}

Recall the decomposition  \eqref{co-martext} and let 
\begin{equation} \label{defmprime}
m' = m +(I-V) g_2 \, . 
\end{equation}
It follows that $(U^i m')_i$ is a stationary  sequence of ${\mathbb L}^1 (\p)$ martingale differences with respect to $(\mathcal F_{i,\infty})_i$.  Since $T$ is ergodic, according to Theorem \ref{EJth}, if 
\begin{equation} \label{condliminfMartP1}
\liminf_{n \rightarrow \infty} \frac{\E (| \sum_{i=0}^{n-1} U^i m' |)}{ \sqrt{n}} < \infty \, ,
\end{equation}
then $m' \in {\mathbb L}^2 (\p)$.  By \eqref{co-martext}, 
\[
\Big \Vert \sum_{i=0}^{n-1} U^i m'  \Big \Vert_1 \leq \Big \Vert \sum_{i=0}^{n-1} U^i f \Big \Vert_1 + 2 \Vert g_1 \Vert_1 + 4   \Vert g_3 \Vert_1 \, .
\]
Hence, since $g_1$ and $g_3$ are in ${\mathbb L}^1(\p)$,
\[
\liminf_{n \rightarrow \infty} \frac{\E (| \sum_{i=0}^{n-1} U^i m' |)}{ \sqrt{n}} \leq \liminf_{n \rightarrow \infty} \frac{\E (| \sum_{i=0}^{n-1} U^i f |)}{ \sqrt{n}} \, , 
\]
which is finite under the first part of  \eqref{condliminf-field-1}. Therefore \eqref{condliminfMartP1} holds and $m' \in {\mathbb L}^2 (\p)$. Next recall that  $m' = m +(I-V) g_2$ and recall that $(V^j m)_j$ is a stationary  sequence of ${\mathbb L}^1 (\p)$ martingale differences with respect to $(\mathcal F_{\infty,j})_j$.  Since $S$ is ergodic, according again to Theorem \ref{EJth}, to prove that $m \in {\mathbb L}^2 (\p)$, it suffices to prove that 
\begin{equation} \label{condliminfMartP2}
\liminf_{n \rightarrow \infty} \frac{\E (| \sum_{j=0}^{n-1} V^j m |)}{ \sqrt{n}} < \infty \, .
\end{equation}
But since $m' = m +(I-V) g_2$ and $g_2$ is in ${\mathbb L}^1(\p)$, proving \eqref{condliminfMartP2} is reduced to show  that 
\begin{equation} \label{condliminfMartP3}
\liminf_{n \rightarrow \infty} \frac{\E (| \sum_{j=0}^{n-1} V^j m' |)}{ \sqrt{n}} < \infty \, .
\end{equation}
With this aim, recall first that $m' \in {\mathbb L}^2 (\p)$ and note that 
\begin{equation} \label{condliminfMartP41}
\frac{1}{\sqrt n} \Big \Vert \sum_{j=0}^{n-1} V^j m'  \Big \Vert_1 \leq \frac{1}{\sqrt n} \Big \Vert \sum_{j=0}^{n-1} V^j m'  \Big \Vert_2 : = \sigma_n \, .
\end{equation}
For any fixed positive integer $n$, let $d := n^{-1/2}  \sum_{j=0}^{n-1} V^j m' $. Since 
$(U^i m')_i$ is a stationary and ergodic  sequence of  martingale differences in ${\mathbb L}^2 (\p)$, so is $(U^i d)_i$.  By the CLT for stationary and ergodic martingales in  ${\mathbb L}^2 (\p)$, as $N \rightarrow \infty$, $ N^{-1/2} \sum_{i=0}^{N-1} U^i d$ converges in distribution to a centered Gaussian random variable $G_n$ with variance $\sigma_n^2$. Hence, by \cite[Th. 3.4]{Bi99} and noticing that $\E |G_n| = \sigma_n \sqrt{2/\pi}$,  for any fixed positive integer $n$, we get 
\begin{equation} \label{condliminfMartP4}
\sigma_n \leq \sqrt{\frac{\pi}{2}}  \liminf_{N \rightarrow \infty}  \frac{1}{\sqrt{N}}  \Big \Vert \sum_{i=0}^{N-1}  U^i  d \Big \Vert_1 = \sqrt{\frac{\pi}{2}}  \liminf_{N \rightarrow \infty}  \frac{1}{\sqrt{n N}}  \Big \Vert \sum_{i=0}^{N-1} \sum_{j=0}^{n-1} U^i V^j m'  \Big \Vert_1 \, .
\end{equation}
But, recalling \eqref{co-martext} and that $ m' = m +(I-V) g_2$, we have 
\begin{multline*}
  \Big \Vert \sum_{i=0}^{N-1} \sum_{j=0}^{n-1} U^i V^j m'  \Big \Vert_1  \leq  \Big \Vert \sum_{i=0}^{N-1} \sum_{j=0}^{n-1} U^i V^j f \Big \Vert_1 
+ \Big \Vert  \sum_{j=0}^{n-1}  (I-U^N)V^j g_1 \Big \Vert_1  \\
+  \Big \Vert  (I-U^N)(I-V^n) g_3 \Big \Vert_1 \, .
\end{multline*}
Hence 
\[
\frac{1}{\sqrt{ N}}  \Big \Vert \sum_{i=0}^{N-1} \sum_{j=0}^{n-1} U^i V^j m'  \Big \Vert_1  \leq \frac{1}{\sqrt{ N}}  \Big \Vert \sum_{i=0}^{N-1} \sum_{j=0}^{n-1} U^i V^j f \Big \Vert_1 
+ \frac{2n }{\sqrt{ N}}  \Vert g_1 \Vert_1  
+  \frac{4}{\sqrt{ N}}   \Vert  g_3  \Vert_1 \, .
\]
Since $g_1$ and $g_3$ are in ${\mathbb L}^1 (\p)$, the two last terms of the right-hand side are converging to zero as $N$ tends to infinity. Hence, taking into account \eqref{condliminfMartP41} and \eqref{condliminfMartP4}, we get
\begin{equation*} \label{condliminfMartP5}
\liminf_{n \rightarrow \infty} \frac{\E (| \sum_{j=0}^{n-1} V^j m' |)}{ \sqrt{n}}  \leq \sqrt{\frac{\pi}{2}}  \liminf_{n \rightarrow \infty}  \liminf_{N \rightarrow \infty}   \frac{1}{\sqrt{nN}}  \Big \Vert \sum_{i=0}^{N-1} \sum_{j=0}^{n-1} U^i V^j f \Big \Vert_1 \, .
\end{equation*}
which is finite by the second part of condition \eqref{condliminf-field-2}.  This ends the proof of \eqref{condliminfMartP3} (and then of \eqref{condliminfMartP2}) and leads to the fact  that $m $ is in ${\mathbb L}^2 (\p)$. Next recall that we have proved that $m'$ defined in \eqref{defmprime} is in ${\mathbb L}^2 (\p)$ which combined with the fact that $m $ is in ${\mathbb L}^2 (\p)$ implies that $(I-V)g_2$ is in ${\mathbb L}^2 (\p)$. On another hand setting 
\[
m'' = m +(I-U) g_1 \, , 
\]
we can use previous arguments to infer that $m''$ is in ${\mathbb L}^2 (\p)$. Hence taking into account that $m  \in {\mathbb L}^2 (\p)$, we get that  $(I-U)g_1$ is in ${\mathbb L}^2 (\p)$. This ends the proof of the theorem.

\subsection{Proof of Theorem \ref{GordinextCE}}

Let $T, S$ be two commuting probability  preserving bijective transformations on $(\Omega, {\mathcal A}, \mu)$, $Uf = f\circ T$, 
$Vf=f\circ S$.
We suppose that ${\mathcal A}$ is generated by $U^iV^je$, $(i,j)\in {\mathbb Z}^2$, and $U^iV^je$ are mutually independent. Denote 
$$
{\mathcal C} = \sigma\{U^i e : i\in {\mathbb Z}\} \, .
$$
$({\mathcal C}, T)$ is thus a Bernoulli dynamical system and the sigma algebras $S^j{\mathcal C}$, $j\in  {\mathbb Z}$, are mutually independent.

Let us recall the so-called Rokhlin lemma. 

\begin{lemma}[Rokhlin lemma]Let $(\Omega, {\mathcal A}, \mu, T)$ be an ergodic dynamical system, $n$ a positive integer,  and 
$\epsilon>0$. Then there exists a set $F\in {\mathcal A}$ such that
\begin{itemize}
\item[(1)] $F, T^{-1}F, \dots, T^{-n+1}F$ are disjoint, \\
\item[(2)] $\mu(\cup_{i=0}^{n-1} T^{-i}F) > 1-\epsilon$.
\end{itemize}
\end{lemma}
The sequence of sets $F, T^{-1}F, \dots, T^{-n+1}F$ defined in the Rokhlin lemma is called  a Rokhlin tower. 

For any integer $k \geq 1$, we set 
\[
  n_k = [2^{k/2} ], \quad m_k = 2^{k} \, .
\]

By the Rokhlin lemma there exists a Rokhlin tower $F, T^{-1}F, \dots, T^{-N_k+1}F$ with $N_k = m_kn_k \sim 2^{3k/2}$. Note that $1 /N_k \geq \mu (T^{-i}F) > (1- \epsilon ) /N_k$, for any $i=0, \ldots, N_k-1$. 
We define
$$
  g_k(\omega) = \begin{cases}
                        (j+1)\sqrt{m_k/n_k}\,\,\,&\text{on}\,\,\,T^{-j}F,\,\,\, j=0\dots,n_k-1, \\
                       (2n_k-j-1)\sqrt{m_k/n_k}\,\,\,&\text{on}\,\,\,T^{-j}F,\,\,\, j=n_k,\dots, 2n_k-1, \\
                       0\,\,\,&\text{on the rest of}\,\,\,\Omega.
\end{cases}
$$
We can notice that 
$$
  g_k - Ug_k = \begin{cases}
                      \sqrt{m_k/n_k}\,\,\,&\text{on}\,\,\,T^{-j}F,\,\,\, j=0\dots,n_k-1, \\
                      -\sqrt{m_k/n_k}\,\,\,&\text{on}\,\,\,T^{-j}F,\,\,\, j=n_k,\dots, 2n_k-1, \\
                       0\,\,\,&\text{on the rest of}\,\,\,\Omega.
\end{cases}
$$
In addition, we have
\begin{align*}
&   \|g_k -Ug_k\|_2^2  \leq  2 \times  \frac{m_k}{n_k} \times  \frac{1}{m_k} \sim  \frac2{2^{k/2}}, \\
 &  \|g_k\|_1  \leq 2 n_k^2  \sqrt{ \frac{m_k}{n_k} }  \times  \frac{1}{n_k m_k} = 2  \sqrt{ \frac{n_k}{m_k} }   \leq  2^{1-k/4}, \\
 &  \Big \Vert  \sum_{i=0}^{2n_k-1} U^i(g_k-Ug_k)  \Big  \Vert_2 =  \Vert g_k - U^{2n_k}g_k \Vert_2  \leq  2 \sqrt{n_k}\, .
  \end{align*}
 Notice that all sums of $U^i(g_k-Ug_k)$ are ${\mathcal C}$-measurable hence the random variables $(V^j(g_k - U^{2n_k}g_k))_j$
are iid. \newline
The support of $g_k - U^{2n_k}g_k$ is included in the union $B_k$ of $F,\dots, T^{-4n_k+1}F$ hence is of measure 
$\leq  4/ m_k = 2^{-k+2}$. \newline
Next, let 
\[
A_k = \{ |g_k - U^{2n_k}g_k|  \geq (1/2)\sqrt{m_k n_k}\}  \, .
\]
It follows that the set $A_k$ is included in $B_k$ and is of measure $2 / m_k= 2^{-k+1}$. 
Because the sigma algebras $S^i {\mathcal C}$, $i  \in {\mathbb Z}$,  are mutually independent, the sets $S^{-j}A_k$, $j=0,\dots,m_k-1$, are independent.
 For $h = 1_{A_k}$ using that ${\rm e}^{2 \ln (1-x) /x } \geq {\rm e}^{-4 }$ for any $x \in ]0, 1/2]$ and that  $m_k  \mu(A_k)=2$, we thus have, for any $k \geq 2$, 
\[
  \mu\Big( \sum_{j=0}^{m_k-1} V^jh = 1 \Big )  = m_k  \mu(A_k) (1- \mu(A_k))^{m_k-1} \geq  2/e^4 \, .
\]
We then  conclude that 
\begin{equation} \label{Conclusion1CE}
 \mu\Big(\frac1{\sqrt{n_k}} \frac1{\sqrt{m_k}} \Big|\sum_{i=0}^{2n_k-1} \sum_{j=0}^{m_k-1} U^iV^j (g_k - Ug_k)\Big|
  \geq 1/2\Big) \geq 2/e^4.
\end{equation}
By recursion we shall define a strictly increasing sequence $k_\ell \nearrow\infty$ and then set 
\[
  g = \sum_{\ell=1}^\infty g_{k_\ell} \text{ and } f = g-Ug \, .
\]
\medskip

1. For $\ell=1$ we put $k_\ell=1$. 
\smallskip

2. Suppose that for $1\leq \ell'<\ell$ the $k_{\ell'}$ have been defined. \newline
All the functions $g_{k_{\ell'}}$ are bounded ($0\leq g_{k_{\ell'}} \leq \sqrt{n_{k_{\ell'}}m_{k_{\ell'}} }$)   hence the sums 
\[
  \sum_{i=0}^{n-1} \sum_{\ell'=1}^{\ell-1} U^i(g_{k_{\ell'}} - Ug_{k_{\ell'}}) =
  \sum_{\ell'=1}^{\ell-1} (g_{k_{\ell'}} - U^ng_{k_{\ell'}}),
\]
$n\geq 1$, are uniformly bounded. \newline
If $k_\ell$ is sufficiently large we thus get
\[
  \frac1{\sqrt{n_{k_\ell}}} \sum_{\ell'=1}^{\ell-1}   \Big|\sum_{i=0}^{2n_{k_\ell}-1}
  U^i(g_{k_{\ell'}} - Ug_{k_{\ell'}})\Big| < 1/2^\ell \, .
\]
Note that, for $j\in {\mathbb Z}$, the functions 
\[
  V^j \Big ( \frac1{\sqrt{n_{k_\ell}}} \sum_{\ell'=1}^{\ell-1} \sum_{i=0}^{2n_{k_\ell}-1}
  U^i(g_{k_{\ell'}} - Ug_{k_{\ell'}}) \Big ) 
\]
are martingale differences. Hence
\begin{equation} \label{Conclusion2CE*}
  \Big\| \frac1{\sqrt{n_{k_\ell}}} \frac1{\sqrt{m_{k_\ell}}} \sum_{i=0}^{2n_{k_\ell}-1} \sum_{j=0}^{m_{k_\ell}-1} 
  U^iV^j \sum_{\ell'=1}^{\ell-1} (g_{k_{\ell'}} - Ug_{k_{\ell'}}) \Big\|_2 \leq \frac1{2^\ell}.
\end{equation}
\medskip
Recall now that  $\|g_k -Ug_k\|_2 \leq  \sqrt2/2^{k/4}$. Hence choosing $k_\ell$ sufficiently large we get
$$
  \Big\| \sum_{i=0}^{2n_{k_{\ell'}}-1} U^i(g_{k_\ell} - Ug_{k_\ell}) \Big\|_2 \leq \frac1{4^\ell}
$$
for all $1\leq \ell'<\ell$.
\medskip

Having constructed the sequence of $k_\ell$ this way we thus have
$$
    \Big\|\frac1{\sqrt{n_{k_\ell}}} \sum_{\ell'=\ell+1}^{\infty}  \sum_{i=0}^{2n_{k_\ell}-1}
  U^i(g_{k_{\ell'}} - Ug_{k_{\ell'}}) \Big\|_2 < 1/2^\ell.
$$
Using the fact that, for $j\in {\mathbb Z}$, the functions 
\[
  V^j \Big ( \frac1{\sqrt{n_{k_\ell}}} \sum_{\ell'=\ell+1}^{\infty} \sum_{i=0}^{2n_{k_\ell}-1}
  U^i(g_{k_{\ell'}} - Ug_{k_{\ell'}}) \Big ) 
\]
are martingale differences, we then derive that  
\begin{equation} \label{Conclusion2CE**}
  \Big\| \frac1{\sqrt{n_{k_\ell}}} \frac1{\sqrt{m_{k_\ell}}} \sum_{i=0}^{2n_{k_\ell}-1} \sum_{j=0}^{m_{k_\ell}-1} 
  U^iV^j \sum_{\ell'=\ell+1}^{\infty} (g_{k_{\ell'}} - Ug_{k_{\ell'}}) \Big\|_2 \leq \frac1{2^\ell} \, .
\end{equation}
Then, the upper bounds  \eqref{Conclusion2CE*} and  \eqref{Conclusion2CE**} entail that 
\begin{equation} \label{Conclusion2CE}
  \Big\| \frac1{\sqrt{n_{k_\ell}}} \frac1{\sqrt{m_{k_\ell}}} \sum_{i=0}^{2n_{k_\ell}-1} \sum_{j=0}^{m_{k_\ell}-1} 
  U^iV^j  (I-U) (g - g_{k_{\ell}}) \Big\|_2 \leq \frac2{2^\ell}.
\end{equation}
Hence taking into account  \eqref{Conclusion1CE} and \eqref{Conclusion2CE}, it follows that, for $f = g - Ug$, the sequence 
$(n_1 n_2)^{-1/2}  \sum_{i=0}^{n_1-1} \sum_{j=0}^{n_2-1} U^iV^j f $ cannot converge in distribution to zero. 

\medskip

Next note that for any $p$ and $q$ fixed, by independence, 
\[
\sum_{i =0}^p  \sum_{j =0}^q {\mathbb E} ( U^iV^j (g - Ug) | {\mathcal F}_{0,0})  =  \sum_{i =0}^p  {\mathbb E}  ( U^i (g - Ug) | {\mathcal F}_{0,0} )
 = g -  {\mathbb E}  ( U^{p+1} g  | {\mathcal F}_{0,0} )  \, .
 \]
But, by the construction of $g$, $\lim_{p \to \infty} \Vert  U^{p+1} g  \Vert_1 =0$. Hence $\Vert  {\mathbb E}  ( U^{p+1} g  | {\mathcal F}_{0,0} ) \Vert_1$  is going to zero as $p \rightarrow \infty$. Therefore condition \eqref{L1-cond-field} is satisfied with $f=g-Ug$. 

\medskip

It remains to prove that the conditions  \eqref{condliminf-field-1} and \eqref{condliminf-field-2} are satisfied with $f=g-Ug$. With this aim, note first that
$$
  \frac1{\sqrt{n}} \sum_{i=0}^{n-1} U^i f  = \frac1{\sqrt{n}} (g-U^ng) \to 0
$$\
in ${\mathbb L}^1$ for $n\to\infty$ (recall that $g\in {\mathbb L}^1$). Next, since the random variables $( V^j f)_{j \geq 0}$ are independent and square integrable,  
\[\frac1{\sqrt{m}} \Big  \Vert  \sum_{j=0}^{m-1} V^j f   \Big \Vert_2 =  \Vert g-U g \Vert_2 < \infty \, .
\] 
Hence both conditions in  \eqref{condliminf-field-1} are satisfied.  On another hand, for every $m$ fixed, $$
  \frac1{\sqrt{n}} \frac1{\sqrt{m}} \Big \Vert \sum_{i=0}^{n-1} \sum_{j=0}^{m-1} U^iV^j (g - Ug)\Big \Vert_1  \leq    \frac{ 2 \sqrt m}{\sqrt{n}}   \Vert g \Vert_1 \to_{n \rightarrow \infty} 0 \, , 
$$
proving the second part of condition  \eqref{condliminf-field-2}. It remains to prove its first part.  Here we use particular properties of $g$ constructed above. 
We have found a sequence of 
$n_{k_\ell}$ for which there exists a positive constant $c>0$ such that 
\[
  \frac1{\sqrt{n_{k_\ell}}} \Big \| \sum_{i=0}^{2n_{k_\ell}-1} U^i(g - Ug) \Big \|_2 \leq c \big ( 2^{- \ell} + 1 \big ) \, .
\]
Since $ \big ( V^j  \sum_{i=0}^{2n_{k_\ell}-1} U^i(g - Ug) \big )_{j \geq 0}$ is a stationary sequence of   martingale differences in ${\mathbb L}^2$, it follows that 
\[
  \frac1{\sqrt{n_{k_\ell}}} \frac1{\sqrt{m}} 
  \Big\|\sum_{i=0}^{2n_{k_\ell}-1} \sum_{j=0}^{m-1} U^iV^j (g - Ug)\Big\|_2 \leq c \big ( 2^{- \ell} + 1 \big ) \leq 2 c  \, .
\]
Because the upper bound is uniform for all $n_{k_\ell}$, the first part of   condition  \eqref{condliminf-field-2} holds true.

\subsection{Proof of Theorem \ref{Gordinextpositive}}

We shall use the coboundary decomposition \eqref{co-martext}.  Note first  that $(U^i V^j m)_{i,j}$ is an ergodic and stationary ${\mathbb L}^2 (\p)$ orthomartingale field. Then, according to the CLT for   ergodic fields of martingale differences as obtained in \cite{V15}, as $\min(n_1,n_2) \rightarrow \infty$,  the sequence
$(n_1n_2)^{-1/2}\sum_{i=0}^{n_1-1} \sum_{j=0}^{n_2-1} U^i V^j m $ converges in distribution  to a centered Gaussian random variable with variance $\Vert m \Vert_2^2$.   Theorem \ref{Gordinextpositive}  then follows from the following proposition. 
 \begin{proposition} \label{Gordinextpositiveprop}  Assume the conditions of Theorem  \ref{Gordinext} and that  condition \eqref{condlimSomme} is satisfied. 
 Then 
 \begin{equation}  \label{consGordinextpositiveprop}
 \lim_{\min(n_1,n_2) \rightarrow \infty } \frac{ \Vert S_{n_1,n_2} (f)  -  S_{n_1,n_2} (m)  \Vert_1}{ \sqrt{n_1 n_2}}  = 0 \, .
\end{equation}
\end{proposition}
{\bf Proof of Proposition \ref{Gordinextpositiveprop}.} Since
\[
(n_1n_2)^{-1/2} \Big \Vert \sum_{i=0}^{n_1-1} \sum_{j=0}^{n_2-1} U^i V^j (I-U)(I-V) g_3 \Big \Vert_1  \leq 4  (n_1n_2)^{-1/2}  \Vert g_3 \Vert_1  \rightarrow_{\min(n_1,n_2) \rightarrow \infty} 0 \, , 
\]
the convergence \eqref{consGordinextpositiveprop}  will follow if one can prove that, as $\min(n_1,n_2) \rightarrow \infty$,
\begin{equation} \label{neglig1}
(n_1n_2)^{-1/2} \Big (  \Big \Vert \sum_{i=0}^{n_1-1} \sum_{j=0}^{n_2-1} U^i V^j (I-U)g_1 \Big \Vert_1 + \Big \Vert \sum_{i=0}^{n_1-1} \sum_{j=0}^{n_2-1} U^i V^j (I-V)g_2 \Big \Vert_1 \Big )  \rightarrow 0 \, .
\end{equation}
Since $(V^j (I-U) g_1)_{j \geq 0} $ and $(U^i (I-U) g_2)_{j \geq 0} $ are sequences of martingale differences in ${\mathbb L}^2 (\p)$, we shall rather prove \eqref{neglig1}  in  ${\mathbb L}^2 (\p)$ and show that 
\begin{equation} \label{neglig1-2}
\lim_{n \rightarrow \infty }\frac{ \Vert (I-U^n) g_1 \Vert_2}{{\sqrt n}}  =0 \,   \text{ and }  \,  \lim_{n \rightarrow \infty }\frac{ \Vert (I-V^n) g_2 \Vert_2}{{\sqrt n}}  =0  \, .
\end{equation}
With this aim, we start by noticing that, for any $n$ fixed, 
\[
d_{1,n} :=  \frac{1}{\sqrt{n}} \sum_{i=0}^{n-1} U^i  \big ( m + (I-U) g_1 \big ) 
\]
is such that $(V^j d_{1,n} )_{j \geq 0 }$ is  an ergodic and stationary sequence of ${\mathbb L}^2 (\p)$ martingale differences with respect to the filtration $(\mathcal F_{\infty,j})_j$. Hence, by the CLT for ergodic and stationary martingales, as $N \to \infty$, 
\[
\frac{1}{\sqrt N}\sum_{j=1}^N  V^j d_{1,n} \rightarrow^{{\mathcal D}}  G_{1,n} \, , 
\]
where $G_{1,n} $ is a centered random Gaussian with standard deviation  $C_n = \Vert d_{1,n} \Vert_2$. Since $\big ( N^{-1/2} \big | \sum_{j=1}^N  V^j d_{1,n} \big |  \big )_{N \geq 1}$ is uniformly integrable, by the convergence of moments theorem (see \cite[Th. 3.5]{Bi99}) we have in particular that 
\[
\lim_{N \to \infty }  \frac{ 1}{\sqrt{N}}\Big \Vert  \sum_{j=1}^N  V^j d_{1,n}  \Big \Vert_1 = \Vert  G_{1,n} \Vert_1 = \sqrt{\frac{2}{\pi}} C_n \, .
\]
But  $\lim_{N \to \infty }  \frac{ 1}{\sqrt{N}}\big \Vert  \sum_{j=0}^{N-1}  V^j d_{1,n}  \big \Vert_1 = \lim_{N \to \infty }  \frac{ 1}{\sqrt{n N}}\big \Vert S_{n,N} (f)  \big \Vert_1$. So, overall,  for any $n$ fixed,
\[
\lim_{N \to \infty }  \frac{ 1}{\sqrt{n N}}\big \Vert S_{n,N} (f)  \big \Vert_1  =  \sqrt{\frac{2}{\pi}} C_n  \, ,
\]
implying by standard arguments that there exists an increasing subsequence $(n_k)$ tending to infinity such that 
\begin{equation} \label{P1positive}
\lim_{k \to \infty }  \Big |  \frac{ 1}{\sqrt{n_k  k }}\big \Vert S_{n_k,k} (f)  \big \Vert_1  -  \sqrt{\frac{2}{\pi}} C_{n_k} \Big |  =0  \, .
\end{equation}
Next, note that 
\[
C^2_n  =   \Vert m \Vert^2_2  +  \frac{2}{n} { \mathbb E} \Big ( (I-U^n) g_1   \sum_{i=0}^{n-1}  U^i m  \Big ) + \frac{ \Vert (I-U^n) g_1 \Vert^2_2}{n }  \, . 
\] 
Notice first that since $(V^j  (I- U^n) g_1)_{j \geq 0}$ is an ergodic and stationary sequence of ${\mathbb L}^2 (\p)$ martingale differences with respect to the filtration $(\mathcal F_{\infty,j})_{j\geq 0}$, we have, by using \cite[Th. 3.4]{Bi99} and arguments used to get  \eqref{condliminfMartP4}, 
\[
\frac{ \Vert (I-U^n) g_1 \Vert_2}{\sqrt{n} } \leq \sqrt{\frac{\pi}{2}}  \liminf_{N \to \infty} \frac{1}{\sqrt{n N}}  \Big \Vert \sum_{i=0}^{n-1} \sum_{j=0}^{N-1} U^i V^j  (I-U)g_1  \Big \Vert_1 \, .
\]
But, according to the coboundary decomposition \eqref{co-martext}, for any $n$ fixed, 
\[
\liminf_{N \to \infty} \frac{1}{\sqrt{n N}}  \Big \Vert \sum_{i=0}^{n-1} \sum_{j=0}^{N-1} U^i V^j  (I-U)g_1  \Big \Vert_1 = \liminf_{N \to \infty} \frac{1}{\sqrt{n N}}  \Big \Vert \sum_{i=0}^{n-1} \sum_{j=0}^{N-1} U^i V^j  (f-m)  \Big \Vert_1 \, .
\]
In addition,
\[
\frac{1}{\sqrt{n N}}  \Big \Vert \sum_{i=0}^{n-1} \sum_{j=0}^{N-1} U^i V^j  (f-m)  \Big \Vert_1 \leq \frac{1}{\sqrt{n N}}  \Big \Vert \sum_{i=0}^{n-1} \sum_{j=0}^{N-1} U^i V^j  f \Big \Vert_1  + \Vert m \Vert_2\, .
\]
So, overall, taking into account  condition \eqref{condlimSomme}, we get 
\begin{equation} \label{P2positiveante}
\kappa:= \limsup_{n \rightarrow \infty}\frac{ \Vert (I-U^n) g_1 \Vert_2}{\sqrt{n} }  \leq \sqrt{\frac{\pi}{2}}  \Big (  \Vert m \Vert_2 + \lim_{n,N \to \infty} \frac{1}{\sqrt{n N}}  \big \Vert S_{n,N} (f)  \big \Vert_1 \Big  ) < \infty \, .
\end{equation}
Now, for any positive real $A$, write
\begin{multline}  \label{P2positiveante-2}
 \Big | C^2_n   -    \Vert m \Vert^2_2  -  \frac{ \Vert (I-U^n) g_1 \Vert^2_2}{n }  \Big | 
   \leq  \frac{2 A }{\sqrt{n}}  \Vert (I- U^n) g_1 \Vert_1  \\ +   2 \frac{ \Vert (I- U^n) g_1 \Vert_2}{\sqrt{n}}  \times  \Big ( \frac{ 1}{n}   \E   \Big   (  \Big |  \sum_{i=1}^n  U^i m  \Big |^2 {\bf 1}_{\{ | \sum_{i=1}^n  U^i m |  > A \sqrt n \}} \Big )^{1/2} \, . 
\end{multline}
Hence, using that $n^{-1/2} \Vert (I- U^n) g_1 \Vert_1 \to_{n \to \infty} 0$ and taking into account  \eqref{P2positiveante} and the fact that $\big ( n^{-1} (\sum_{i=0}^{n-1}  U^i m )^2 \big )_{n \geq 1}$ is  uniformly integrable, we derive that the terms in the right-hand side of \eqref{P2positiveante-2} tend to zero by letting first $n$ goes to infinity and after $A$. Therefore
\[
\lim_{n \rightarrow \infty} \Big | C^2_n   -    \Vert m \Vert^2_2  -  \frac{ \Vert (I-U^n) g_1 \Vert^2_2}{n }  \Big | = 0 \, .
\]
Assume now that 
\begin{equation} \label{P2positive}
\kappa= \limsup_{n \to \infty }   \frac{ 1}{\sqrt{n }} \Vert (I-U^n) g_1 \Vert_2  >0  \, .
\end{equation}
Then, there exists an increasing subsequence $(n'_\ell)_{\ell \geq 1}$ tending to infinity such that 
\begin{equation} \label{P3positive}
\lim_{\ell \to \infty }    \frac{ 1}{\sqrt{n'_\ell }} \Vert (I-U^{n'_\ell}) g_1 \Vert_2 = \kappa  \text{ and then } \lim_{\ell \to \infty } C^2_{n'_\ell}  =    \Vert m \Vert^2_2 + \kappa^2  \, .
\end{equation}
According to \eqref{P1positive} and \eqref{P3positive}, we then infer that if \eqref{P2positive} holds true then there exist two increasing subsequences $(n^{\prime \prime}_\ell)_{\ell \geq 1}$ and $(k^{\prime\prime}_\ell)_{\ell \geq 1}$ tending to infinity such that 
\begin{equation} \label{P4positive}
\lim_{\ell \to \infty }  \frac{ 1}{\sqrt{n^{\prime\prime}_\ell k^{\prime \prime}_\ell }}\Big \Vert S_{n^{\prime \prime}_\ell , k^{\prime \prime}_\ell} (f)  \Big \Vert_1   >   \sqrt{\frac{2}{\pi}} \Vert m \Vert_2  \, .
\end{equation}
But, using once again  the coboundary decomposition \eqref{co-martext},  note that 
\[
\frac{1}{n} S_{n,n} (f) = \frac{1}{n} S_{n,n} (m)  +  \frac{1}{n} (I-U^n )  \sum_{j=1}^n V^j g_1 +   \frac{1}{n} (I-V^n ) \sum_{i=1}^n U^i g_2 +  \frac{1}{n} (I-V^n ) (I-V^n ) g_3 \, .
\]
Using that $g_3$ is in ${\mathbb L}^1 (\p )$, and the fact that since $g_1$ and $g_2$ are  in ${\mathbb L}^1 (\p )$, the Birkhoff theorem in ${\mathbb L}^1 (\p )$ implies that 
\[
\lim_{n\to \infty }  \frac{ 1}{ n }\Vert(I-U^n )  \sum_{j=1}^n V^j g_1 \Vert_1 =0 \text{ and } \lim_{n\to \infty }  \frac{ 1}{ n }\Vert(I-V^n ) \sum_{i=1}^n U^i g_2 \Vert_1 =0 \, , 
\]
we get 
\begin{equation} \label{P5positive}
\lim_{n\to \infty }  \frac{ 1}{ n }\Vert S_{n,n} (f) \Vert_1  = \lim_{n\to \infty }  \frac{ 1}{ n }\Vert S_{n,n} (m) \Vert_1  \, .
\end{equation}
But since $n^{-1}  S_{n,n} (m)$ converges in distribution  to a centered Gaussian random variable with variance $\Vert m \Vert_2^2$ and $(\frac{ 1}{ n }\vert S_{n,n} (m) \vert)_{n \geq 1}$ is  uniformly integrable, we derive, by the convergence of moments theorem,  that 
\[
\lim_{n\to \infty }  \frac{ 1}{ n }\Vert S_{n,n} (m) \Vert_1 =  \sqrt{\frac{2}{\pi}}  \Vert m \Vert_2 \, .
\]
This result together with  \eqref{P5positive}  imply that 
\begin{equation} \label{P6positive}
\lim_{n\to \infty }  \frac{ 1}{ n }\Vert S_{n,n} (f) \Vert_1  = \sqrt{\frac{2}{\pi}}  \Vert m \Vert_2  \, .
\end{equation}
Clearly, under condition  \eqref{condlimSomme},  if \eqref{P2positive} is supposed to be true, \eqref{P4positive} and \eqref{P6positive} are not compatible. This proves that \eqref{P2positive} cannot be true and then that the first part of \eqref{neglig1-2} is satisfied. With similar arguments, one can prove that, provided the additional condition  \eqref{condlimSomme} is assumed,  the second part of \eqref{neglig1-2} is also satisfied.  This ends the proof of the proposition and then of the theorem.

\section{Extension to multidimensional index of higher dimension} \label{sectiongeneral}

To state the extension of Theorems  \ref{Gordinext} and  \ref{Gordinextpositive} to higher dimensions, some additional notations are needed. Let $d \geq 1$ and $(T_{\underline{i}})_{\underline{i} \in {\mathbb Z}^d}$ be $ {\mathbb Z}^d$ actions on  $(\Omega, {\mathcal A}, \p)$ generated by commuting invertible and measure-preserving transformations $T_{\varepsilon_q}$, $1 \leq q \leq d$. Here ${\varepsilon}_q$ is the vector of $ {\mathbb Z}^d$ which has $1$ ate the $q$-th place and $0$ elsewhere. By $U_{\underline{i}}$  we denote the operator in ${\mathbb L}^p$ ($1 \leq p \leq \infty$) defined by $ U_{\underline{i}} f = f \circ T_{\underline{i}}$, $\underline{i} \in {\mathbb Z}^d$.  By $\underline{i} \leq \underline{j}$, we understand $i_k \leq j_k$ for all $1 \leq k \leq d$. 

We suppose that there is a completely commuting $({\mathcal  F}_{\underline{j}})_{\underline{j} \in {\mathbb Z}^d}$ , i.e. there is a $\sigma$-algebra ${\mathcal F}_{\underline{0}}$ such that ${\mathcal  F}_{\underline{i}}=T_{-\underline{i}} {\mathcal F}_{\underline{0}}$, for $\underline{i} \leq \underline{j}$ we have ${\mathcal  F}_{\underline{i}} \subset {\mathcal  F}_{\underline{j}}$ and for an integrable $f$,
\[
 \E ( \E ( f | {\mathcal  F}_{i_1, \ldots, i_d }) | {\mathcal  F}_{j_1, \ldots, j_d } )  =  \E ( f | {\mathcal  F}_{i_1 \wedge j_1, \ldots, i_d  \wedge j_d}) \, .
\]

By ${\mathcal F}_{\ell}^{(k)}$ we denote the $\sigma$-algebra generated by all ${\mathcal  F}_{\underline{i}}$ with $\underline{i} = (i_1, \ldots, i_d)$ with $i_k \leq \ell $ and $i_j \in {\mathbb Z}$ for $1 \leq j \leq d$, $j \neq \ell$.   For $\sigma$-algebras $ {\mathcal G} \subset  {\mathcal F}  \subset  {\mathcal A} $, by ${\mathbb L}^p ( \mathcal F) \ominus {\mathbb L}^p ( \mathcal G)$ we denote the space of $f  \in   {\mathbb L}^p ( \mathcal F) $ for which $\E ( f | {\mathcal G}) =0$ a.s. 

For  $f$  a $\ {\mathcal F}_{\underline{0}}$-measurable centered ${\mathbb L}^1(\p)$ random variable,  it has been proved in Voln\'y \cite[Th. 4]{V18} that the condition
\begin{equation} \label{L1-cond-fieldgeneral}
\text{ the series } \sum_{ i_1, \ldots, i_d \geq 0} \E ( U_{i_1, \ldots, i_d} f | {\mathcal F}_{\underline{ 0}}) \text{ converges in ${\mathbb L}^1(\p)$}
\end{equation}
ensures the existence of the following  orthomartingale-coboundary decomposition:  
\begin{equation} \label{OMCdecgenneral}
f = m + \sum_{\emptyset \subsetneq  J  \subsetneq <d> }  \prod_{s \in J}(I-U_{s})m_J + \prod_{s=1}^d(I-U_{s}) g 
\end{equation}
where $<d> := \{1, \ldots , d\}$, $m$, $g$ and $m_J$ belong to ${\mathbb L}^1({\mathcal F}_{\underline{0}},  \p)$,  ${\mathbb L}^1( \prod_{s=1}^d T_s {\mathcal F}_{\underline{0}},  \p)$ and ${\mathbb L}^1( \prod_{s \in J} T_s {\mathcal F}_{\underline{0}},  \p)$ respectively and $(U^{{\underline i}}_{< d>} m  )_{\underline{i} \in {\mathbb Z}^d}$  and $(U^{{\underline i}}_{J^c} m_J  )_{\underline{i} \in {\mathbb Z}^{d-|J|}}$ are  orthomartingale diffferences random fields  for $\emptyset \subsetneq  J  \subsetneq <d> $.  

\medskip

For any positive integer $k$, define ${\mathcal S}_k$ the set of all the permutations of $\{1, \ldots, d\}$. We are now in position to state the extension of Theorems  \ref{Gordinext} and  \ref{Gordinextpositive} .

\begin{theorem} \label{Gordinextgeneral} Let  $d \geq 1$ and $f$  a $\ {\mathcal F}_{\underline{0}}$-measurable centered ${\mathbb L}^1(\p)$ random variable. Let  ${\underline n}_d = (n_1, \ldots, n_d)$ and  $S_{\underline{n}_d} (f) = \sum_{i_1=0}^{n_1-1}  \ldots   \sum_{i_d=0}^{n_d-1}
U_{i_1, \ldots, i_d}f$. 
Assume that  each of the  transformations $T_{\varepsilon_q}$, $1 \leq q \leq d$, is ergodic. Suppose also that   condition \eqref{L1-cond-fieldgeneral}  holds and that for any integer $k \in \{1, \ldots, d\}$ and all  the permutations $\sigma$ in ${\mathcal S}_k$, 
\begin{equation} \label{extensioncondliminf-field} 
\liminf_{ n_{\sigma (1)} \to \infty} \cdots  \liminf_{ n_{\sigma (k)} \to \infty}  \frac{\E \big ( \big | \sum_{i_1=0}^{n_{\sigma (1)}-1}  \ldots   \sum_{i_k=0}^{n_{\sigma (k)}-1}U_{i_1, \ldots, i_k}f \big  | \big )}{ (\prod_{i=1}^k n_{\sigma (i)})^{1/2} } < \infty \, .
\end{equation}
Then $m \in {\mathbb L}^2 ( \p)$ and for any set $J$ such that $ \emptyset \subsetneq  J  \subsetneq <d>$, $ \prod_{s \in J}(I-U_{s})m_J \in {\mathbb L}^2 ( \p)$ ($m$ and $m_J$ are  defined in \eqref{OMCdecgenneral}).  
If, in addition, \begin{equation} \label{condlimfieldTCL}
\lim_{\min(n_1,n_2, \ldots, n_d) \rightarrow \infty} \frac{\E (|S_{{\underline n}_d } (f)  |)}{ (\prod_{i=1}^d n_i)^{1/2} } \text{ exists,} 
\end{equation}
then $(\prod_{i=1}^d n_i)^{-1/2} S_{{\underline n}_d } (f) $ converges in  distribution to a centered normal variable (that can be degenerate) as $\min(n_1,n_2, \ldots, n_d) \rightarrow \infty$. \end{theorem}

\noindent {\bf Proof of Theorem \ref{Gordinextgeneral}.} The result will follow by recurrence. Note that it holds for $d=1$ and also for $d=2$ as shown in the previous section. Assume that it holds for $d-1$ and let us prove it for $d$. Recall the decomposition  \eqref{OMCdecgenneral} and let 
\begin{equation} \label{defmprime}
m'=  m + \sum_{\emptyset \subsetneq  J  \subseteq <d>_1 }  \prod_{s \in J}(I-U_{s})m_J  \, , 
\end{equation}
where $<d>_1 = <d> \backslash \{1\} =   \{2, \ldots , d\}$. Note that $(U_1^i m')_{i \in {\mathbb Z}} $ is a stationary sequence of ${\mathbb L}^1 (\p)$ martingale differences w.r.t. $({\mathcal F}_{i, 0, \ldots,0})_{i \in {\mathbb Z}}$.   Since $T_1$ is ergodic, according to Theorem \ref{EJth}, if 
\[
\liminf_{n \rightarrow \infty} \frac{\E (| \sum_{i=0}^{n-1} U_1^i m' |)}{ \sqrt{n}} < \infty \, ,
\]
then $m' \in {\mathbb L}^2 (\p)$.  This follows from  the decomposition \eqref{OMCdecgenneral} and the fact that, by condition \eqref{extensioncondliminf-field}, $\liminf_{n \rightarrow \infty}\frac{\E (| \sum_{i=0}^{n-1} U_1^i f |)}{ \sqrt{n}} < \infty $.   Next, starting from \eqref{defmprime}, and taking into account  the induction hypothesis, namely: Theorem \ref{Gordinextgeneral} holds for $d-1$, we infer that if  for any  integer $k \in \{2, \ldots, d\}$ and all  the permutations $\sigma$ in ${\mathcal S}_k$,   
\begin{equation} \label{P1general} 
\liminf_{ n_{\sigma (1)} \to \infty} \cdots  \liminf_{ n_{\sigma (k)} \to \infty}  \frac{\E \big ( \big | \sum_{i_1=0}^{n_{\sigma (1)}-1}  \ldots   \sum_{i_k=0}^{n_{\sigma (k)}-1}U_{i_1, \ldots, i_k} m' \big  | \big )}{ (\prod_{i=1}^k n_{\sigma (i)})^{1/2} } < \infty \, .
\end{equation}
then $m \in {\mathbb L}^2 ( \p)$  and, for any set $J$ such that $ \emptyset \subsetneq  J  \subsetneq <d>_1$, $ \prod_{s \in J}(I-U_{s})m_J \in {\mathbb L}^2 ( \p)$.  By using similar arguments as those developed in the proof of Theorem \ref{Gordinext}, we infer that \eqref{P1general} is satisfied under condition \eqref{extensioncondliminf-field}. Hence  $m \in {\mathbb L}^2 ( \p)$. Then, using in addition that  $m' \in {\mathbb L}^2 ( \p)$, we conclude that, for any set $J$ such that $ \emptyset \subsetneq  J  \subseteq <d>_1$, $ \prod_{s \in J}(I-U_{s})m_J \in {\mathbb L}^2 ( \p)$. The first part of Theorem \ref{Gordinextgeneral}  follows by using $d-1$ times the same arguments and replacing  $ <d>_1$ by $ <d>_i$ for $i=2, \ldots, k$. The second part of the theorem follows by applying  
 the CLT for ergodic and stationary fields of  orthomartingales as proved in Voln\'y  \cite{V15} for $(\prod_{i=1}^d n_i)^{-1/2} S_{{\underline n}_d } (m) $ and by using similar arguments as those developed in the proof of Theorem \ref{Gordinextpositive}.


\begin{thebibliography}{99}     

 
\bibitem{Bi99} Billingsley, P. 
{\it Convergence of probability measures.}
Second edition. Wiley Series in Probability and Statistics: Probability and Statistics.  John Wiley $\&$ Sons, Inc., New York, 1999. x+277 pp.

\bibitem{D98} Dedecker, J. A central limit theorem for stationary random fields. \textit{Probab.  Theory  Related Fields} 110 (1998), no. 3,  397--426.

\bibitem{EG} El Machkouri, M. and Giraudo, D.  Orthomartingale-coboundary decomposition for stationary random fields.  \textit{Stoch. Dyn.} 16 (2016), no. 5, 1650017, 28 pp.

\bibitem{Gi18} Giraudo, D.   Invariance principle via orthomartingale approximation.  \textit{Stoch. Dyn.} 18 (2018), no. 6, 1850043, 29 pp.
 
\bibitem{Go73} Gordin, M. I. (1973). Abstracts of Communication, T.1:A-K, International Conference on Probability Theory, Vilnius.

\bibitem{EJ-85}  Esseen, C.G. and  Janson, S. 
On moment conditions for normed sums of independent variables and martingale differences.
\textit{Stochastic Process. Appl}. 19 (1985), no. 1, 173--182.

\bibitem{PZ17} Peligrad, M. and Zhang, N. On the normal approximation for random fields via martingale methods. \textit{Stochastic Process. Appl}. 128 (2018), no. 4, 1333--1346.

\bibitem{V15}  Voln\'y, D.  A central limit theorem for fields of martingale differences. \textit{C. R. Math. Acad. Sci. Paris} 353 (2015), no. 12, 1159--1163. 

\bibitem{V18} Voln\'y, D. Martingale-coboundary representation for stationary random fields. \textit{Stoch. Dyn.} 18 (2018), no. 2, 1850011, 18 pp.

\bibitem{V13} Voln\'y, D.
Approximating martingales and the central limit theorem for strictly stationary processes. \textit{Stochastic Process. Appl.} 44 (1993), no. 1, 41--74. 

\bibitem{VW14} Voln\'y, D. and Wang, Y.  An invariance principle for stationary
random fields under Hannan's condition. \textit{Stochastic Process. Appl.}  124 (2014), no. 12,  4012--4029.

\bibitem{WV13}  Wang, Y. and Woodroofe, M. A new condition for the invariance principle for stationary random fields. \textit{Statist. Sinica }23 (2013), no. 4, 1673--1696.

\end{thebibliography}
\end{document}